\documentclass[11pt]{article}

\def\R{{\bf R}}
\def\N{{\bf N}}

\def\dy{\displaystyle}
\def\ep{{\varepsilon}}

\def\-{|\!|\!|}

\newtheorem{thm}{Theorem}
\newtheorem{lm}{Lemma}[section]
\newtheorem{prop}{Proposition}[section]

\newtheorem{re}{Remark}[section]
\makeatletter
 
 \@addtoreset{equation}{section}
\makeatother

\usepackage{latexsym,graphicx,}%showkeys}
\usepackage[usenames]{color}

\title{The blow-up and lifespan of solutions to systems of semilinear wave equation
with critical exponents in high dimensions}
\author{Yuki Kurokawa $^a$ \quad Hiroyuki Takamura $^b$ \quad Kyouhei Wakasa $^c$}
\date{\begin{center}
{\small {\it 
$^a$ General Education, 
Yonago National College of Technology, 
4448 Yonago, Tottori 683-8502, Japan. 
e-mail : kurokawa@yonago-k.ac.jp\\
$^b$ Department of Complex and Intelligent Systems,
Faculty of Systems Information Science, 
Future University Hakodate, 
116-2 Kamedanakano-cho, Hakodate, Hokkaido 041-8655, Japan. 
e-mail : takamura@fun.ac.jp\\
$^c$ The 1st year of Graduate School of Systems Information Science, 
Future University Hakodate, 
116-2 Kamedanakano-cho, Hakodate, Hokkaido 041-8655, Japan. 
e-mail : g2111045@fun.ac.jp}}
\end{center}}

\begin{document}
\maketitle
\begin{abstract}

\par
In this paper we prove the blow-up theorem in the critical case
for weakly coupled systems of semilinear wave equations in high dimensions. 
The upper bound of the lifespan of the solution is precisely clarified. 
\end{abstract}
\begin{center}
\begin{tabular}{l}
\end{tabular}
\end{center}
%%%%%%%%%%%%%%%%%%%%%%%%%%
\section{Introduction}
\label{sec:Intro}
\par
Let us consider the following systems of semilinear wave equations; 
\begin{equation}
\left\{
\begin{array}{l}
u_{tt}-\Delta u=|v|^p\\
v_{tt}-\Delta v=|u|^q
\end{array}
\right.
\qquad\mbox{in}\quad \R^n\times[0,\infty)
\label{system}
\end{equation}
for $p,q>1$ and $n\ge 2$ with the data of the form;
\begin{equation}
\left\{
\begin{array}{l}
u(x,0)=\ep f_1(x),\ u_t(x,0)=\ep g_1(x)\\
v(x,0)=\ep f_2(x),\ v_t(x,0)=\ep g_2(x)
\end{array}
\right.
\qquad\mbox{for}\quad x\in\R^n,
\label{data}
\end{equation}
where $\ep>0$ is a small parameter.
We assume that $f_1,f_2,g_1,g_2\in C_0^\infty(\R^n)$ for the simplicity. 
\par
In order to describe the results on (\ref{system}), we set
\begin{equation}
\label{F}
F(p,q,n)\equiv\max\left\{\dy\frac{q+2+p^{-1}}{pq-1},
\ \dy\frac{p+2+q^{-1}}{pq-1}\right\}-\dy\frac{n-1}{2}.
\end{equation}
DelSanto, Georgiev and Mitidieri \cite{DGM} first showed
that the system (\ref{system}) with (\ref{data}) for $n\ge 2$
has a global in time solution for sufficiently small $\ep$ if $F(p,q,n)<0$,
while a solution for some positive data blows up in finite time if $F(p,q,n)>0$. 
When the blow-up occurs, it is known that the maximal time $T(\ep)$,
so-called \lq\lq lifespan", of the existence of solutions for arbitrarily fixed data
can be estimated as 
\begin{equation}
\label{lifespan:sub}
c\ep^{-F(p,q,n)^{-1}}\le T(\ep)\le C\ep^{-F(p,q,n)^{-1}}, 
\end{equation}
where $c$ and $C$ are positive constants independent of $\ep$.
See \cite{AKT, GTZ, Kuro1, Kuro2, KuTa}. 
\par
When $F(p,q,n)=0$, the non-existence of global in time solutions,
namely $T(\ep)<\infty$, was shown
by DelSanto and Mitidieri \cite{DM} for $n=3$. 
Moreover, a sharp estimate of the lifespan;
\begin{equation}
\label{lifespan:critical}
\begin{array}{cl}
\exp\left(c\ep^{-\min\{p(pq-1),q(pq-1)\}}\right)\le T(\ep)
\le \exp\left(C\ep^{-\min\{p(pq-1),q(pq-1)\}}\right)\quad
&\mbox{for}\quad p\neq q,\\
\exp\left(c\ep^{-p(p-1)}\right)\le T(\ep)\le \exp\left(C\ep^{-p(p-1)}\right)\quad
&\mbox{for}\quad p= q,\\
\end{array}
\end{equation}
was obtained for $n=2,3$ by \cite{AKT, KO1, KuTa}. 
In high dimensions, $n\ge 4$, only the lower bounds of the lifespan are estimated.
See \cite{GTZ, Kuro1, Kuro2}. 
The upper bounds of the lifespan have not been obtained,
because we have similar technical difficulties to the corresponding problem
for single equations;
\begin{equation}
\label{single}
\left\{
\begin{array}{l}
u_{tt}-\Delta u=|u|^p\qquad\mbox{in}\quad \R^n\times [0,\infty)\\
u(x,0)=\ep f(x),\ u_t(x,0)=\ep g(x)
\end{array}
\right.
\end{equation}
for $n\ge 4$ and $p=p_0(n)$.
To show the blow-up result in this case was an open problem for a long period. 
See \cite{DG, GLS, G2, G3, J1, LZ, Li, Sc, Si2, St, Z1, Z2}, etc.
for the other cases on (\ref{single}). 
Here, $p_0(n)$ is a positive root of the quadratic equation $2+(n+1)p-(n-1)p^2=0$. 
Note that $F(p,p,n)=0$ is equivalent to $p=p_0(n)$. 
\par
Finally, Yordanov and Zhang \cite{YZ} or Zhou \cite{Z3}
independently obtained a blow-up result of $T(\ep)<\infty$ for this open case. 
Later, Takamura and Wakasa \cite{TW} have succeeded to get
the optimal estimate of the lifespan
by introducing a new iteration argument based on the method in \cite{YZ}.
Zhou Yi and Han Wei \cite{ZH} have recently reproved the theorem in \cite{TW}
along with the method in \cite{Z3}. 
\par 
Our aim in this article is to clarify the lifespan of a solution
to the system (\ref{system}) for the critical case 
in high space dimensions by employing the argument in \cite{TW}.  
\begin{thm}
\label{thm:lifespan}
Let $n\ge 4$ and $F(p,q,n)=0$ with $1<p\le q$. 
Assume that
\begin{equation}
\label{supp:ID}
\begin{array}{l}
(f_1,g_1)\in H^1(\R^n)\cap L^q(\R^n)\times L^2(\R^n),
\ (f_2,g_2)\in H^1(\R^n)\times L^2(\R^n),\\
f_1(x)=f_2(x)=g_1(x)=g_2(x)=0\ \mbox{for}\ |x|>R>0
\end{array}
\end{equation}
and that $f_1,f_2,g_1,g_2$ are non-negative, especially
that $g_1$ and $g_2$ do not vanish identically.
Moreover, suppose that the problem (\ref{system}) with (\ref{data}) has a solution
$(u,v)\in C([0,T(\ep));H^1(\R^n)\cap L^q(\R^n)\times H^1(\R^n)$ with 
$(u_t,v_t)\in(C([0,T(\ep));L^2(\R^n)))^2$ satisfying 
\begin{equation}
\label{support}
\mbox{supp}(u,v,u_t,v_t)\subset\{(x,t)\in\R^n\times[0,T(\ep))\ ;\ |x|\le t+R\}. 
\end{equation}
Then, there exists a positive constant $\ep_0=\ep_0(f_1,f_2,g_1,g_2,p,q,n,R)$
such that $T(\ep)$ has to satisfy 
\begin{equation}
\label{lifespan}
T(\ep)\le\left\{
\begin{array}{ll}
\exp(C\ep^{-p(pq-1)})\ &\mbox{if}\quad p\neq q\\
\exp(C\ep^{-p(p-1)})\ &\mbox{if}\quad p=q,
\end{array}
\right.
\end{equation}
for $0<\ep\le \ep_0$, where $C$ is a positive constant independent of $\ep$. 
\end{thm}
\begin{re}
When $n\ge 4$ and $F(p,q,n)=0$ with $1<p\le q$, we have that $p\le 2$. 
Because $p_0(4)=2$ and $p_0(n)$ is monotonously decreasing in $n$. 
Moreover, we can take $R\ge 1$ in the proof of Theorem \ref{thm:lifespan}
without loss of the generality. 
\end{re}
\begin{re}
The counter case, $1<q\le p$, can be obtained by symmetricity of (\ref{system}).
In fact, if one substitutes $p$ with $q$ and $u$ with $v$,
the upper bounds in (\ref{lifespan:critical}) immediately follow
from Theorem\ref{thm:lifespan}. 
\end{re}
%%%%%%%%%%%%%%%%%%%%%%%%%%%%%%%%%%%%%%
\par
The proof is based on the argument in \cite{TW}. 
A blow-up property of a nonlinear system of the ordinary differential inequalities
will be employed to show a blow-up of a solution to (\ref{system}). 
The growing up of $L^p$ norm of $v$ is crucial 
to get the sharp estimate of the lifespan. 
In order to obtain the growth of the norm by iteration argument,
we will employ the integral inequalities of $L^p$ norm of $v$
and $L^q$ norm of $u$ as a frame of the iteration,
which are based on the method of Yordanov and Zhang \cite{YZ} .
Such an argument was introduced by Takamura and Wakasa \cite{TW}.  
\par
This paper is organized as follows. 
In the next section, we shall show a blow-up property for nonlinear systems
of ordinary differential inequalities. 
By making use of this, we prove Theorem \ref{thm:lifespan} in Section \ref{section:proof}. 
In Section \ref{section:iteration},
we complete the iteration argument. 
%%%%%%%%%%%%%%%%%%%%%%%%%%%%%%%%%%%
\section{Blow-up for systems of ODIs with a critical balance}
\label{section:ODI}
\par
As stated in Introduction,
we shall show a blow-up theorem for ordinary differential inequalities
with a critical balance in exponents. 
\begin{lm}
\label{lm:ODIs}
Let $p,q>1$, $a>0$, $\alpha,\beta\ge 0$ and $\alpha+p\beta=a(pq-1)+2(p+1)$. 
Suppose that $U,V\in C^2([0,T])$ satisfy 
\begin{eqnarray}
&&U(t)\ge Kt^a \hspace{1cm} \mbox{for}\quad t\ge T_0,\label{ineq1}\\
&&U''(t)\ge A(t+R)^{-\alpha}|V(t)|^p \hspace{1cm} \mbox{for}\quad t\ge 0,\label{ineq2}\\
&&V''(t)\ge B(t+R)^{-\beta}|U(t)|^q \hspace{1cm} \mbox{for}\quad t\ge 0,\label{ineq3}\\
&&U(0)\ge0,\ U'(0)>0,\ V(0)\ge0,\ V'(0)>0,\label{ineq4}
\end{eqnarray}
where all $A,B,K,R,T_0$ are positive constants with $T_0\ge R$. 
Then, $T$ must satisfy that $T\le 2T_1$ provided $K\ge K_0$, where 
\begin{equation}
K_0=\left[
\dy\frac{a^{2p+2}2^{\alpha+p\beta+2p+1}(q+2)^{2p+2}(q+1)^p}{AB^p}
\left(1-\dy\frac{1}{2^{a\delta}}\right)^{-2p-2}
\right]^{\frac{1}{pq-1}}
\label{K0}
\end{equation}
and
\begin{equation}
T_1=\max\left\{T_0,\ \dy\frac{(2q+3)U(0)}{U'(0)}\right\}
\label{T1}
\end{equation}
with a positive constant $\delta\in(0,(pq-1)/(2p+2))$. 
\end{lm}
%%%%%%%%%%%%%%%%%%%%%%%%%%%%%%%%%5
{\it Proof.}\ 
We shall prove this lemma by contradiction. 
Assume that $T>2T_1$.
We note that 
\begin{equation}
U(t)\ge0,\ U'(t)>0,\ V(t)\ge0\ \mbox{and}\ V'(t)>0\quad\mbox{for}\quad t\ge 0 
\label{ineq5}
\end{equation}
by (\ref{ineq2}), (\ref{ineq3}) and (\ref{ineq4}). 
Multiplying (\ref{ineq3}) by $U'(t)$ and integrating it over $[0,t]$, we have 
\[
\begin{array}{l}
U'(t)V'(t)-U'(0)V'(0)-\dy\int_0^tU''(s)V'(s)ds\\
\ge B\dy\int_0^t(s+R)^{-\beta}U(s)^qU'(s)ds\\
\ge B(t+R)^{-\beta}\dy\int_0^tU(s)^qU'(s)ds\\
= \dy\frac{B}{q+1}(t+R)^{-\beta}\{U(t)^{q+1}-U(0)^{q+1}\}.
\end{array}
\]
Then, (\ref{ineq5}) gives us 
\[
U'(t)V'(t)>\dy\frac{B}{q+1}(t+R)^{-\beta}\{U(t)^{q+1}-U(0)^{q+1}\}. 
\]
Multiplying this inequality by $U'(t)$ again and integrating it over $[0,t]$, we have 
\[
\begin{array}{l}
U'(t)^2V(t)-U'(0)^2V(0)-2\dy\int_0^tU'(s)U''(s)V(s)ds\\
>\dy\frac{B}{q+1}\dy\int_0^t(s+R)^{-\beta}\{U(s)^{q+1}-U(0)^{q+1}\}U'(s)ds\\
\ge\dy\frac{B}{q+1}(t+R)^{-\beta}\dy\int_0^t\{U(s)^{q+1}-U(0)^{q+1}\}U'(s)ds\\
=\dy\frac{B}{q+1}(t+R)^{-\beta}\phi(t), 
\end{array}
\]
where we set
\begin{equation}
\phi(t)=\dy\frac{U(t)^{q+2}-U(0)^{q+2}}{q+2}-U(0)^{q+1}\{U(t)-U(0)\}
\ge 0\qquad\mbox{for}\quad t\ge 0.
\label{phi}
\end{equation}
Thus, by (\ref{ineq2}) and (\ref{ineq5}), we have
\[
V(t)>\dy\frac{B}{q+1}(t+R)^{-\beta}\phi(t)U'(t)^{-2}.
\]
Substituting this expression into (\ref{ineq2}), we obtain that 
\[
U''(t)U'(t)^{2p}>A\left(\dy\frac{B}{q+1}\right)^p
(t+R)^{-\alpha-p\beta}\phi(t)^p\qquad\mbox{for}\ t\ge 0.
\]
Multiplying this inequality by $U'(t)$ and integrating it over $[0,t]$, we have 
\begin{equation}
\begin{array}{l}
\label{est1}
\dy\frac{1}{2p+2}\{U'(t)^{2p+2}-U'(0)^{2p+2}\}\\
>A\left(\dy\frac{B}{q+1}\right)^p\dy\int_0^t
(s+R)^{-\alpha-p\beta}\phi(s)^pU'(s)ds\nonumber\\
\ge A\left(\dy\frac{B}{q+1}\right)^p
(t+R)^{-\alpha-p\beta}\dy\int_0^t\phi(s)^p\phi'(s)\{\phi'(s)\}^{-1}U'(s)ds.
\end{array}
\end{equation}
\par
By (\ref{phi}), one can easily see that 
\[
\{\phi'(s)\}^{-1}U'(s)=\{U(s)^{q+1}-U(0)^{q+1}\}^{-1}
\ge U(t)^{-q-1}\qquad\mbox{for}\quad 0<s\le t. 
\]
Then, $s$-integral in (\ref{est1}) is estimated from below as follows. 
\begin{equation}
s\mbox{-integral}
\ge U(t)^{-q-1}\dy\int_0^t\phi(s)^p\phi'(s)ds=
\dy\frac{1}{p+1}U(t)^{-q-1}\phi(t)^{p+1}.
\label{est2}
\end{equation}
Moreover, by the monotonicity of $U(t)$, we have 
\begin{eqnarray*}
\phi(t)
&\ge& \dy\frac{U(t)^{q+1}\{U(t)-U(0)\}}{q+2}-U(0)^{q+1}\{U(t)-U(0)\}\\
&=& \{U(t)-U(0)\}\left\{\dy\frac{U(t)^{q+1}}{q+2}-U(0)^{q+1}\right\}\\
&\ge& \{U(t)-U(0)\}U(t)^q\left\{\dy\frac{U(t)}{q+2}-U(0)\right\}. 
\end{eqnarray*}
It also follows from the monotonicity of $U'(t)$ that $U(t)-U(0)\ge U'(0)t$.
Thus we have 
\[
U(t)-U(0)
\ge \dy\frac{U(t)}{q+2}-U(0)
\ge \dy\frac{U(t)}{2(q+2)}+\dy\frac{U'(0)t+U(0)}{2(q+2)}-U(0)
\ge \dy\frac{U(t)}{2(q+2)}
\]
provided $t\ge(2q+3)U(0)/U'(0)$. Hence we get
\begin{equation}
\phi(t)\ge\dy\frac{1}{4(q+2)^2}U(t)^{q+2}
\qquad\mbox{for}\quad t\ge \dy\frac{(2q+3)U(0)}{U'(0)}. 
\label{est3}
\end{equation}
\par
Therefore, it follows from (\ref{est1}), (\ref{est2}) and (\ref{est3}) that 
\[
U'(t)>\left[\dy\frac{2A}{\{4(q+2)^2\}^{p+1}}
\left(\dy\frac{B}{q+1}\right)^p\right]^{\frac{1}{2p+2}}
(t+R)^{-\frac{\alpha+p\beta}{2p+2}}U(t)^{\frac{pq+2p+1}{2p+2}}
\quad\mbox{for}\ t\ge\frac{(2q+3)U(0)}{U'(0)}.   
\]
By the definition of $K_0$ in (\ref{K0}), this inequality can be rewritten as
\[
U'(t)>aK_0^{-\frac{pq-1}{2p+2}}
\left(1-\dy\frac{1}{2^{a\delta}}\right)^{-1}
t^{-\frac{\alpha+p\beta}{2p+2}}U(t)^{\frac{pq+2p+1}{2p+2}}
\]
if $t\ge R$. 
Now, we assume that $t\ge T_1$, where $T_1$ is the one in (\ref{T1}). 
Multiplying the last inequality by $U(t)^{-1-\delta}>0$
with a constant $\delta\in(0,(pq-1)/(2p+2))$ and replacing $U(t)$
in the right hand side by (\ref{ineq1}), we have 
\[
U(t)^{-1-\delta}U'(t)
>a\left(\dy\frac{K}{K_0}\right)^{\frac{pq-1}{2p+2}}
K^{-\delta}\left(1-\dy\frac{1}{2^{a\delta}}\right)^{-1}t^{-1-a\delta}
\]
because of the critical balance, $\alpha+p\beta=a(pq-1)+2(p+1)$. 
Integrating the above inequality over $[T_1,t]$, we have 
\[
\dy\frac{1}{\delta}\left\{U(T_1)^{-\delta}-U(t)^{-\delta}\right\}
>\dy\frac{1}{\delta}\left(\dy\frac{K}{K_0}\right)^{\frac{pq-1}{2p+2}}
K^{-\delta}\left(1-\dy\frac{1}{2^{a\delta}}\right)^{-1}
\left(T_1^{-a\delta}-t^{-a\delta}\right). 
\]
By the assumption of $T>2T_1$, one can set $t=2T_1$. 
Then, neglecting the second term in the left hand side, we get 
\[
1>\left(\dy\frac{K}{K_0}\right)^{\frac{pq-1}{2p+2}}
\left\{\dy\frac{U(T_1)}{KT_1^a}\right\}^\delta. 
\] 
By making use of (\ref{ineq1}) again for the right hand side,
we find that this inequality implies 
\[
1>\left(\dy\frac{K}{K_0}\right)^{\frac{pq-1}{2p+2}}.
\]
This contradicts to $K\ge K_0$. 
Therefore we have $T\le 2T_1$.
Lemma \ref{lm:ODIs} is now established. 
\hfill$\Box$
%%%%%%%%%%%%%%%%%%%%%%%%%%%%%%%%%%%
\section{Proof of Theorem \ref{thm:lifespan}}
\label{section:proof}
\par
In this section we shall prove Theorem \ref{thm:lifespan} by using Lemma \ref{lm:ODIs}. 
Let us define 
\begin{equation}
U(t)=\dy\int_{\R^n}u(x,t)dx\quad\mbox{and}\quad 
V(t)=\dy\int_{\R^n}v(x,t)dx,
\label{UV}
\end{equation}
where $(u,v)$ is a solution to (\ref{system}) with (\ref{data}) satisfying (\ref{support}). 
\par 
First, we shall show (\ref{ineq2})-(\ref{ineq4}) in Lemma \ref{lm:ODIs} for $U$ and $V$ in (\ref{UV}). 
Integrating two equations of (\ref{system}) in $x\in\R^n$, we have 
\[
U''(t)=\dy\int_{\R^n}|v(x,t)|^p\ dx\quad\mbox{and}\quad 
V''(t)=\dy\int_{\R^n}|u(x,t)|^q\ dx
\]
by the support property (\ref{support}). 
By making use of H\"{o}lder's inequality together
with the support property (\ref{support}) again, we have 
\[
\left\{
\begin{array}{l}
U''(t)\ge B_n^{1-p}(t+R)^{-n(p-1)}|V(t)|^p\\
V''(t)\ge B_n^{1-q}(t+R)^{-n(q-1)}|U(t)|^q
\end{array}
\right.
\quad\mbox{for}\quad t\ge 0, 
\]
where $B_n$ stands for the volume of the unit ball in $\R^n$. 
This implies that (\ref{ineq2}) and (\ref{ineq3}) are valid with
\begin{equation}
\alpha=n(p-1),\ \beta=n(q-1),\ A=B_n^{1-p},\ B=B_n^{1-q}.
\label{parameter}
\end{equation}
The assumption on the positiveness of data (\ref{supp:ID}) gives us (\ref{ineq4}) with
\[
\begin{array}{l}
U(0)=\ep\dy\int_{\R^n}f_1(x)dx\ge0,\qquad U'(0)=\ep\dy\int_{\R^n}g_1(x)dx>0,\\
V(0)=\ep\dy\int_{\R^n}f_2(x)dx\ge0,\qquad V'(0)=\ep\dy\int_{\R^n}g_2(x)dx>0.\\
\end{array}
\]
\par
Next, we shall show the inequality (\ref{ineq1}) in this situation
employing the following estimates of $U''(t)$ up to the case of $p\neq q$, or $p=q$. 
\begin{prop}
\label{prop:iteration1}
Suppose that the assumptions in Theorem \ref{thm:lifespan} are fulfilled and $p<q$. 
Then, $U(t)$ satisfies the following inequality. 
\begin{equation}
\label{iteration1}
U''(t)\ge C_j(t-a_jR)^{n-1-(n-1)p/2}
\left(\log\dy\frac{t+(a_j^2-2)R}{(a_j+2)(a_j-1)R}\right)^{\frac{(pq)^j-1}{pq-1}}
\end{equation}
for $t\ge a_jR$ and $j=1,2,3\cdots$. 
Here we set $a_j=6\cdot 4^{j-1}-2$ and 
\begin{equation}
\label{Cj}
C_1=\tilde{C}\ep^{p^2q},\quad
C_j=\exp\{(pq)^{j-1}(\log C_1+S_j)\}\ \mbox{for}\ j\ge 2, 
\end{equation}
where $\tilde{C}=\tilde{C}(p,q,n)$ is a positive constant
and $S_j$ is a convergent sequence independent of $\ep$ and $t$. 
\end{prop}
\begin{prop}
\label{prop:iteration2}
Suppose that the assumptions in Theorem \ref{thm:lifespan} are fulfilled and $p=q$. 
Then, $U(t)$ satisfies the following inequality. 
\begin{equation}
\label{iteration2}
U''(t)\ge D_j(t-b_jR)^{n-1-(n-1)p/2}
\left(\log\dy\frac{t+(b_j^2-6)R}{(b_j-2)(b_j+3)R}\right)^{\frac{p^{2j}-1}{p-1}}
\end{equation}
for $t\ge b_jR$ and $j=1,2,3\cdots$. 
Here we set $b_j=10\cdot 4^{j-1}-2$ and 
\begin{equation}
\label{Dj}
D_1=\tilde{C}\ep^{p^3},\quad
D_j=\exp\{p^{2(j-1)}(\log D_1+S_j')\}\ \mbox{for}\ j\ge 2,
\end{equation}
where $\tilde{C}=\tilde{C}(p,q,n)$ is a positive constant
and $S_j'$ is a convergent sequence independent of $\ep$ and $t$. 
\end{prop}
\par
These propositions are proved in the next section.
From now on, we shall apply (\ref{iteration1}) and (\ref{iteration2})
to the proof of (\ref{ineq1}). 
We are concentrated on the case $p<q$ only
since the proof for the case of $p=q$ is analogue. 
%%%%%%%%%%%%%%%%%%%%%%%%%%%
\par
First, let $t\ge a_jR$. 
Integrating (\ref{iteration1}) over $[a_jR,t]$, we have 
\begin{equation}
U'(t)-U'(a_jR)\ge C_j\dy\int_{a_jR}^t(s-a_jR)^{n-1-(n-1)p/2}
\left(\log\dy\frac{s+(a_j^2-2)R}{(a_j+2)(a_j-1)R}\right)^{\frac{(pq)^j-1}{pq-1}}ds.
\label{U'}
\end{equation}
We can neglect the second term in the left hand side by (\ref{ineq5}). 
Restricting the time interval to $t\ge (a_j+1)R$,
we can replace the lower limit by $a_jt/(a_j+1)$ because of
\[
a_j R\le \dy\frac{a_j}{a_j+1}t<t\qquad\mbox{for}\quad t\ge (a_j+1)R.
\]
Then it follows that 
\[
s-a_jR\ge \dy\frac{a_j}{a_j+1}t-a_jR=\frac{a_j}{a_j+1}\{t-(a_j+1)R\}
\ge \frac{t-(a_j+1)R}{2}
\]
and
\[
\frac{s+(a_j^2-2)R}{(a_j+2)(a_j-1)R}
\ge \dy\frac{\frac{a_j}{a_j+1}t+(a_j^2-2)R}{(a_j+2)(a_j-1)R}
=\frac{a_jt+(a_j+1)(a_j^2-2)R}{(a_j+1)(a_j+2)(a_j-1)R}
\ge \frac{t+(a_j+1)^2R}{(a_j+1)(a_j+2)R}
\]
for $a_jt/(a_j+1)\le s\le t$. 
Thus, we get by (\ref{U'}) that
\begin{eqnarray*}
U'(t)
&\ge& 
C_j\left(\dy\frac{t-(a_j+1)R}{2}\right)^{n-1-(n-1)p/2}
\left(\log\dy\frac{t+(a_j+1)^2R}{(a_j+1)(a_j+2)R}\right)^{\frac{(pq)^j-1}{pq-1}}\dy\int_{a_jt/(a_j+1)}^tds\\
&\ge& 
\dy\frac{C_j}{2^{n-(n-1)p/2}a_j}\{t-(a_j+1)R\}^{n-(n-1)p/2}
\left(\log\dy\frac{t+(a_j+1)^2R}{(a_j+1)(a_j+2)R}\right)^{\frac{(pq)^j-1}{pq-1}}
\end{eqnarray*}
for $t\ge(a_j+1)R$ because of $a_j=6\cdot 4^{j-1}-2\ge 4$. 
Integrating the last inequality over $[(a_j+1)R,t]$
and treating it in the similar way as above, we have 
\[
U(t)\ge\dy\frac{C_j}{2^{2n+1-(n-1)p}a_j^2}\{t-(a_j+2)R\}^{n+1-(n-1)p/2}
\left(\log\dy\frac{t+(a_j+1)(a_j+2)R}{(a_j+2)^2R}\right)^{\frac{(pq)^j-1}{pq-1}} 
\]
for $t\ge (a_j+2)R$. 
\par
Let us again restrict the time interval to $t\ge (a_j+2)^4R^2$. 
Then it holds that 
\[
\frac{t+(a_j+1)(a_j+2)R}{(a_j+2)^2R}
\ge\sqrt{t}\times\frac{\sqrt{t}}{(a_j+2)^2R}\ge\sqrt{t}
\]
and
\[
t-(a_j+2)R=\frac{t}{2}+\frac{t-2(a_j+2)R}{2}\ge\frac{t}{2}. 
\]
It follows from these inequalities and $a_j=6\cdot 4^{j-1}-2\le 2\cdot 4^j$ that 
\begin{equation}
\label{U}
U(t)\ge K_j(t)t^{n+1-(n-1)p/2}\qquad\mbox{for}\quad t\ge(a_j+2)^4R^2,
\end{equation}
where we set
\begin{equation}
\label{Kj}
K_j(t)=\dy\frac{C_jA}{16^j}
\left(\log\sqrt{t}\right)^{\frac{(pq)^j-1}{pq-1}},
\quad
A=2^{-3n-4+3(n-1)p/2}.
\end{equation}
\par
Now, we consider (\ref{Kj}) for $t\in[(a_j+2)^4R^2,\ (a_{j+1}+2)^4R^2]$. 
Then, by the definition of $C_j$ (\ref{Cj}), $K_j$ in (\ref{Kj}) can be rewritten as 
\[
K_j(t)=\exp\left\{(pq)^{j-1}\log L_j(t)
+\log A-j\log 16-\dy\frac{1}{pq-1}\log(\log\sqrt{t})\right\},
\]
where 
\[
L_j(t)=C_1e^{S_j}(\log\sqrt{t})^{\frac{pq}{pq-1}}. 
\]
Since $S_j$ converges to a certain number,
there exists a constant $S=S(p,q,n)$ such that $S_j\ge S$ for any $j=1,2,3,\cdots$. 
It follows from the definition of $C_1$, (\ref{Cj}), that $L_j\ge e$ holds provided 
\begin{equation}
\label{lifespan1}
\ep^{p(pq-1)}\log t\ge E, 
\end{equation}
where $E=2(\tilde{C}^{-1}e^{1-S})^{(pq-1)/pq}$. 
\par
We assume (\ref{lifespan1}). 
Then it follows from $t\le (a_{j+1}+2)^4R^2$ that 
\[
K_j(t)\ge \exp\left\{(pq)^{j-1}+\log A-j\log 16-\dy\frac{1}{pq-1}\log(\log((a_{j+1}+2)^2R))\right\}. 
\]
Hence we can see that $K_j(t)$ goes to infinity if $j$ tends to infinity. 
Therefore, for $K_0$ defined in (\ref{K0}) with (\ref{parameter}), $a=n+1-(n-1)p/2$ 
and a constant $\delta\in(0,(pq-1)/(2p+2))$, 
there exists an integer $J=J(f_1,f_2,g_1,g_2,n,p,q,R)$ such that 
\[
K_j(t)\ge K_0\qquad\mbox{for}\quad t\in[(a_j+2)^4R^2,(a_{j+1}+2)^4R^2],
\]
as far as $j\ge J$. 
This implies that 
\[
U(t)\ge K_0t^{n+1-(n-1)p/2}\qquad\mbox{for}\quad t\ge(a_J+2)^4R^2, 
\]
provided (\ref{lifespan1}) is valid. 
\par
Now, we are in a position to prove Theorem\ref{thm:lifespan}
by making use of Lemma\ref{lm:ODIs}. 
Set 
\begin{equation}
\label{T0}
T_0(\ep)=\exp(E\ep^{-p(pq-1)}),
\end{equation}
where $E$ is the one in (\ref{lifespan1}). 
Then there exists a positive constant $\ep_0=\ep_0(f_1,f_2,g_1,g_2,p,q,n,R)$ such that 
\begin{equation}
T_0(\ep)\ge(a_J+2)^4R^2\quad\mbox{and}\quad
2\max\left\{T_0(\ep),\ \dy\frac{(2q+3)U(0)}{U'(0)}\right\}\le\exp(2E\ep^{-p(pq-1)})
\label{ep}
\end{equation}
holds for $0<\ep\le \ep_0$. 
As we see,
(\ref{ineq1}) is now established for $t\ge T_0(\ep)$ with this $\ep$. 
We also obtain other inequalities in Lemma\ref{lm:ODIs}
with (\ref{parameter}) and $a=n+1-(n-1)p/2$. 
Note that the condition in Lemma\ref{lm:ODIs} $\alpha+p\beta=a(pq-1)+2(p+1)$ is
equivalent to the critical relation $F(p,q,n)=0$.
In this way, when $T(\ep)>T_0(\ep)$, (\ref{ineq1}) holds for $t\in[T_0(\ep),T(\ep))$. 
Hence Lemme\ref{lm:ODIs} and (\ref{ep}) show that 
\[
t\le 2\max\left\{T_0(\ep),\ \dy\frac{(2q+3)U(0)}{U'(0)}\right\}\le \exp(2E\ep^{-p(pq-1)}). 
\]
Taking a supremum over $t\in[T_0(\ep),T(\ep))$, we get 
\[
T(\ep)\le\exp(2E\ep^{-p(pq-1)})\qquad\mbox{for}\quad 0<\ep\le\ep_0. 
\]
When $T(\ep)\le T_0(\ep)$, (\ref{lifespan}) is trivial.
Therefore the proof of Theorem \ref{thm:lifespan} is ended.
%%%%%%%%%%%%%%%%%%%%%%%%%%%%%%%
\section{Iteration argument}
\label{section:iteration}
\par
In this section, we will prove
Proposition \ref{prop:iteration1} and Proposition\ref{prop:iteration2} by iteration argument. 
The following integral expressions of $U''$ and $V''$ are the frame in our iteration. 
\begin{prop}
\label{prop:frame}
Suppose that the assumptions in Theorem \ref{thm:lifespan} are fulfilled. 
Then, $U(t)$ and $V(t)$ satisfy
\begin{equation}
U''(t)\ge C\dy\int_0^{t-R}\frac{\rho^{n-1-(n-1)p/2}d\rho}{(t-\rho +R)^{(n-1)p/2}}
\left(\dy\int_0^{(t-\rho-R)/2}V''(s)ds\right)^p,
\label{frame_u}
\end{equation}
\begin{equation}
V''(t)\ge \dy\frac{C}{(t+R)^{(n-1)q/2}}\dy\int_0^{t-R}
\frac{\rho^{n-1}d\rho}{(t-\rho +R)^{(n-1)q/2}}
\left(\dy\int_0^{(t-\rho-R)/2}U''(s)ds\right)^q
\label{frame_v}
\end{equation}
for $t\ge R$, where $C$ is a positive constant independent of $\ep$ and $t$. 
\end{prop}
{\it Proof.} 
Recall that $p\le 2$. 
Then, these inequalities can be immediately obtained by (2.14) and (2.21)
in Yordanov and Zhang \cite{YZ}. 
Therefore we shall omit the proof.
\par
The first step of the iteration is the following estimate. 
\begin{prop}
Suppose that the assumptions in Theorem \ref{thm:lifespan} are fulfilled. 
Then, there exists a positive constant $C=C(f_2,g_2,n,p,R)$ such that 
\begin{equation}
\label{first}
U''(t)\ge C\ep^p(t+R)^{n-1-(n-1)p/2}\qquad\mbox{for}\quad t\ge 0.
\end{equation}
\end{prop}
{\it Proof.} 
This inequality can be proved by the same way as (2.5') in Yordanov and Zhang \cite{YZ}. 
The key estimates, (2.4) and Lemma 2.2 in \cite{YZ}, are obatined
by the first and second equations in (\ref{system}).\\
%%%%%%%%%%%%%%%%%%%%%%%%%%%%%%%%%%%%
\par
Let us continue to prove Proposition\ref{prop:iteration1} by 
making use of the two propositions above. 
Substituting (\ref{first}) into $U''(s)$ in the $s$-integral in (\ref{frame_v}), we have 
\[
V''(t)\ge \dy\frac{C^{q+1}\ep^{pq}}{(t+R)^{(n-1)q/2}}
\dy\int_0^{t-R}\frac{\rho^{n-1}d\rho}{(t-\rho +R)^{(n-1)q/2}}
\left(\dy\int_0^{(t-\rho-R)/2}(s+R)^{n-1-(n-1)p/2}ds\right)^q.
\]
for $t\ge R$. 
Putting the upper limit of the $s$-integral into a part of the negative power of $s+R$,
we have
\[
\begin{array}{ll}
s\mbox{-integral}
&\dy\ge\left(\frac{t-\rho+R}{2}\right)^{-(n-1)p/2}\int_0^{(t-\rho-R)/2}s^{n-1}ds\\
&\dy\ge\frac{(t-\rho-R)^n}{2^nn(t-\rho+R)^{(n-1)p/2}}.
\end{array}
\]
Hence we get
\[ 
\begin{array}{ll}
V''(t)
&\ge\dy\frac{C^{q+1}\ep^{pq}}{2^{nq}n^q(t+R)^{(n-1)q/2}}
\dy\int_0^{t-R}\frac{\rho^{n-1}(t-\rho-R)^{nq}}{(t-\rho +R)^{(n-1)q(p+1)/2}} d\rho\\
&\ge\dy\frac{C^{q+1}\ep^{pq}}{2^{nq}n^q(t+R)^{(n-1)q(p+2)/2}}
\dy\int_0^{t-R}\rho^{n-1}(t-\rho-R)^{nq} d\rho.
\end{array}
\]
Cutting the domain of the $\rho$-integral, we have 
\[
\rho\mbox{-integral}
\ge\left(\dy\frac{t-R}{2}\right)^{nq}\dy\int_0^{(t-R)/2}\rho^{n-1} d\rho
=\dy\frac{(t-R)^{nq+n}}{2^{nq+n}n}. 
\]
Thus we get 
\begin{equation}
\label{v1}
V''(t)\ge \dy\frac{C'(t-R)^{nq+n}}{(t+R)^{(n-1)q(p+2)/2}}
\qquad\mbox{for}\quad t\ge R,
\end{equation}
where $C'=C^{q+1}\ep^{pq}2^{-2nq-n}n^{-q-1}$. 
Next we shall set $t\ge 3R$ and substitute (\ref{v1}) into $V''(s)$ in (\ref{frame_u}). 
Then, it follows that
\[
\begin{array}{ll}
U''(t)
&\ge CC'^p\dy\int_0^{t-3R}\frac{\rho^{n-1-(n-1)p/2}d\rho}{(t-\rho +R)^{(n-1)p/2}}
\left(\dy\int_R^{(t-\rho-R)/2}\dy\frac{(s-R)^{nq+n}}{(s+R)^{(n-1)q(p+2)/2}}ds\right)^p\\
&\ge CC'^p\dy\int_0^{t-3R}\frac{\rho^{n-1-(n-1)p/2}d\rho}{(t-\rho +R)^{(n-1)p(pq+2q+1)/2}}
\left(\dy\int_R^{(t-\rho-R)/2}(s-R)^{nq+n}ds\right)^p\\
&=\dy\frac{CC'^p}{(nq+n+1)^p2^{npq+(n+1)p}}
\dy\int_0^{t-3R}\frac{\rho^{n-1-(n-1)p/2}(t-\rho-3R)^{npq+(n+1)p}}{(t-\rho +R)^{(n-1)p(pq+2q+1)/2}}d\rho
\end{array}
\]
for $t\ge 3R$. 
Here we again restrict the time interval to $t\ge a_1R=4R$. 
Then, it follows from 
\[
5(t-\rho-3R)\ge t-\rho+R\quad\mbox{for}\quad \rho\le t-4R
\]
and
\[
(n-1)p(pq+2q+1)/2-npq-(n+1)p=1\quad\mbox{for}\quad F(p,q,n)=0
\]
that
\[
\begin{array}{ll}
\rho\mbox{-integral}
&\ge\dy\frac{1}{5^{(n-1)p(pq+2q+1)/2}}
\dy\int_{(t-4R)/2}^{t-4R}\frac{\rho^{n-1-(n-1)p/2}}{t-\rho-3R}d\rho\\
&\ge\dy\frac{1}{5^{(n-1)p(pq+2q+1)/2}}
(t-4R)^{-(n-1)p/2}\left(\dy\frac{t-4R}{2}\right)^{n-1}
\dy\int_{(t-4R)/2}^{t-4R}\frac{d\rho}{t-\rho-3R}\\
&=\dy\frac{(t-4R)^{n-1-(n-1)p/2}}{5^{(n-1)p(pq+2q+1)/2}2^{n-1}}
\log\dy\frac{t-2R}{2R}.
\end{array}
\]
Since 
\[
\dy\frac{t-2R}{2R}=\dy\frac{t+8t-18R}{18R}\ge \dy\frac{t+32R-18R}{18R}=\dy\frac{t+14R}{18R}
\]
holds for $t\ge 4R$, we get 
\[
U''(t)
\ge\dy C_1(t-4R)^{n-1-(n-1)p/2}\log\dy\frac{t+14R}{18R}, 
\]
where
\[
\begin{array}{ll}
C_1
&=\dy\frac{CC'^p}{(nq+n+1)^p2^{npq+(n+1)p+n-1}5^{(n-1)p(pq+2q+1)/2}}\\
&=\dy\frac{C^{pq+p+1}\ep^{p^2q}}{(nq+n+1)^pn^{p(q+1)}2^{3npq+(2n+1)p+n-1}5^{(n-1)p(pq+2q+1)/2}}
\end{array}
\]
Therefore (\ref{iteration1}) is true for $j=1$. 
%%%%%%%%%%%%%%%%%%%%%%%%%%%%%%
\par
Next we shall show (\ref{iteration1}) by induction. 
Assume that (\ref{iteration1}) for $t\ge a_j R$ holds and $C_j$ is unknown here except for $j=1$  but will be determined later on. 
When $t\ge(2a_j+1)R$, substituting (\ref{iteration1}) into (\ref{frame_v}), we have 
\[
\begin{array}{ll}
V''(t)
&\ge\dy\frac{CC_j^q}{(t+R)^{(n-1)q/2}}
\int_0^{t-(2a_j+1)R}\frac{\rho^{n-1}d\rho}{(t-\rho+R)^{(n-1)q/2}}\\
&\quad\dy\times
\left\{\int_{a_jR}^{(t-\rho-R)/2}(s-a_jR)^{n-1-(n-1)p/2}
\left(\log\frac{s+(a_j^2-2)R}{(a_j-1)(a_j+2)R}\right)^{\frac{(pq)^j-1}{pq-1}}ds\right\}^q\\
&\ge\dy\frac{CC_j^q}{(t+R)^{(n-1)q/2}}
\int_0^{t-(2a_j+1)R}\frac{\rho^{n-1}I_j(\rho,t)^q}{(t-\rho+R)^{(n-1)q(p+1)/2}}d\rho, 
\end{array}
\]
where 
\[
I_j(\rho,t)=\dy\int_{a_jR}^{(t-\rho-R)/2}(s-a_jR)^{n-1}
\left(\log\frac{s+(a_j^2-2)R}{(a_j-1)(a_j+2)R}\right)^{\frac{(pq)^j-1}{pq-1}}ds. 
\]
Here we restricted the time interval to $t\ge(2a_j+2)R$
and cut the domain of $\rho$-integral to be $[0,t-(2a_j+2)R]$. 
Then, it follows from $(2a_j+1)R\le t-\rho-R$ and
\[
a_jR\le\dy\frac{a_j}{2a_j+1}(t-\rho-R)<\dy\frac{t-\rho-R}{2} 
\]
that one can cut the domain of $s$-integral as $[a_j(t-\rho-R)/(2a_j+1), (t-\rho-R)/2]$. 
In this interval, we have
\[
\begin{array}{ll}
\dy\frac{s+(a_j^2-2)R}{(a_j-1)(a_j+2)R}
&\ge\dy\frac{\frac{a_j}{2a_j+1}(t-\rho-R)+(a_j^2-2)R}{(a_j-1)(a_j+2)R}\\
&=\dy\frac{a_j(t-\rho-R)+(2a_j+1)(a_j^2-2)R}{(2a_j+1)(a_j-1)(a_j+2)R}\\
&\ge\dy\frac{(a_j-1)(t-\rho-R)+(2a_j+1)R+(2a_j+1)(a_j^2-2)R}{(2a_j+1)(a_j-1)(a_j+2)R}\\
&=\dy\frac{t-\rho+a_j(2a_j+3)R}{(2a_j+1)(a_j+2)R}
\end{array}
\] 
and
\[
s-a_jR
\ge \dy\frac{a_j}{2a_j+1}(t-\rho-R)-a_jR
\ge \dy\frac{t-\rho-(2a_j+2)R}{3}.
\]
Hence we have
\[
\begin{array}{ll}
I_j(\rho,t)
&\ge\left(\dy\frac{t-\rho-(2a_j+2)R}{3}\right)^{n-1}
\left(\log\dy\frac{t-\rho+a_j(2a_j+3)R}{(2a_j+1)(a_j+2)R}\right)^{\frac{(pq)^j-1}{pq-1}}\\
&\quad\times\left(\dy\frac{1}{2}-\frac{a_j}{2a_j+1}\right)(t-\rho-R)\\
&\ge\dy\frac{1}{2\cdot 3^{n}a_j}(t-\rho-(2a_j+2)R)^n
\left(\log\dy\frac{t-\rho+a_j(2a_j+3)R}{(2a_j+1)(a_j+2)R}\right)^{\frac{(pq)^j-1}{pq-1}}
\end{array}
\]
because of $a_j\ge 1$ for any $j\in \N$. 
Therefore we obtain
\[
\begin{array}{ll}
V''(t)
&\ge\dy\frac{CC_j^q}{2^q\cdot 3^{nq}a_j^q(t+R)^{(n-1)q/2}}\dy\int_0^{t-(2a_j+2)R}
\frac{\rho^{n-1}\{t-\rho-(2a_j+2)R\}^{nq}}{(t-\rho+R)^{(n-1)q(p+1)/2}}\\
&\hspace{150pt}\times\left(\log\dy
\frac{t-\rho+a_j(2a_j+3)R}{(2a_j+1)(a_j+2)R}\right)^{\frac{q((pq)^j-1)}{pq-1}}d\rho\\
&\ge\dy\frac{CC_j^q}{2^q\cdot 3^{nq}a_j^q(t+R)^{(n-1)q(p+2)/2}}\\
&\quad\times\dy\int_0^{(t-(2a_j+2)R)/2}\!\!\!\rho^{n-1}\{t-\rho-(2a_j+2)R\}^{nq}
\left(\log\dy\frac{t-\rho+a_j(2a_j+3)R}{(2a_j+1)(a_j+2)R}\right)^{\frac{q((pq)^j-1)}{pq-1}}d\rho\\
&\ge\dy\frac{CC_j^q}{2^{(n+1)q}\cdot 3^{nq}a_j^q}\cdot
\dy\frac{\{t-(2a_j+2)R\}^{nq}}{(t+R)^{(n-1)q(p+2)/2}}
\left(\log\dy\frac{t+(4a_j^2+8a_j+2)R}{2(2a_j+1)(a_j+2)R}\right)^{\frac{q((pq)^j-1)}{pq-1}}\\
&\quad\times\dy\int_0^{(t-(2a_j+2)R)/2}\rho^{n-1}d\rho\\
&=\dy\frac{C_j'\{t-(2a_j+2)R\}^{nq+n}}{(t+R)^{(n-1)q(p+2)/2}}
\left(\log\dy\frac{t+(4a_j^2+8a_j+2)R}{2(2a_j+1)(a_j+2)R}\right)^{\frac{q((pq)^j-1)}{pq-1}}, 
\end{array}
\]
where 
\begin{equation}
C_j'=\dy\frac{CC_j^q}{n2^{(n+1)q+n}\cdot 3^{nq}a_j^q}.
\label{Cj'}
\end{equation}
\par
When $t\ge (4a_j+5)R$, replacing the $V''(s)$ in the right hand side
in (\ref{frame_u}) by the last inequality above, we have 
\[
\begin{array}{ll}
U''(t)
&\ge CC_j'^p\dy\int_0^{t-(4a_j+5)R}
\frac{\rho^{n-1-(n-1)p/2}}{(t-\rho +R)^{(n-1)p/2}}d\rho\\
&\quad\times 
\left(\dy\int_{(2a_j+2)R}^{(t-\rho-R)/2}
\dy\frac{\{s-(2a_j+2)R\}^{nq+n}}{(s+R)^{(n-1)q(p+2)/2}}
\left(\log\dy\frac{s+(4a_j^2+8a_j+2)R}{2(2a_j+1)(a_j+2)R}\right)^{\frac{q((pq)^j-1)}{pq-1}}
ds\right)^p\\
&\ge CC_j'^p
\dy\int_0^{t-(4a_j+5)R}\frac{\rho^{n-1-(n-1)p/2}J_j(\rho,t)^p}
{(t-\rho +R)^{(n-1)p(pq+2q+1)/2}}d\rho,
\end{array}
\]
where 
\[
J_j(\rho,t)=
\dy\int_{(2a_j+2)R}^{(t-\rho-R)/2}
\{s-(2a_j+2)R\}^{nq+n}
\left(\log\frac{s+(4a_j^2+8a_j+2)R}{2(2a_j+1)(a_j+2)R}\right)^{\frac{q((pq)^j-1)}{pq-1}}ds
\]
Let us again restrict the time interval to $t\ge (4a_j+6)R$
and cut the domain of $\rho$-integral to be $[0,t-(4a_j+6)R]$. 
Then, it follows from $(4a_j+5)R\le t-\rho-R$ and
\[
(2a_j+2)R\le\dy\frac{2a_j+2}{4a_j+5}(t-\rho-R)\le\dy\frac{t-\rho-R}{2}
\]
that one can cut the domain of the $s$-integral
of $J_j(\rho,t)$ to be $[(2a_j+2)(t-\rho-R)/(4a_j+5), (t-\rho-R)/2]$. 
Then, the inequality
\[
\begin{array}{l}
\dy\frac{s+(4a_j^2+8a_j+2)R}{2(2a_j+1)(a_j+2)R}
\ge\dy\frac{\frac{2a_j+2}{4a_j+5}(t-\rho-R)+(4a_j^2+8a_j+2)R}{2(2a_j+1)(a_j+2)R}\\
=\dy\frac{(2a_j+2)(t-\rho-R)+(4a_j+5)(4a_j^2+8a_j+2)R}{2(4a_j+5)(2a_j+1)(a_j+2)R}\\
\ge\dy\frac{(2a_j+1)(t-\rho-R)+(4a_j+5)R+(4a_j+5)(4a_j^2+8a_j+2)R}
{2(4a_j+5)(2a_j+1)(a_j+2)R}\\
=\dy\frac{t-\rho+(8a_j^2+22a_j+14)R}{2(4a_j+5)(a_j+2)R}
\end{array}
\]
holds for $(2a_j+2)(t-\rho-R)/(4a_j+5)\le s\le(t-\rho-R)/2$. 
Thus we have 
\[
\begin{array}{ll}
J_j(\rho,t)
&\ge\dy\left(\log\frac{t-\rho+(8a_j^2+22a_j+14)R}
{2(4a_j+5)(a_j+2)R}\right)^{\frac{q((pq)^j-1)}{pq-1}}\\
&\quad\dy\times\int_{\frac{2a_j+2}{4a_j+5}(t-\rho-R)}^{(t-\rho-R)/2}
\left\{s-\frac{2a_j+2}{4a_j+5}(t-\rho-R)\right\}^{nq+n}ds\\
&\ge\dy\frac{(t-\rho-R)^{nq+n+1}}{(nq+n+1)(12a_j)^{nq+n+1}}
\left(\log\dy\frac{t-\rho+(8a_j^2+22a_j+14)R}
{2(4a_j+5)(a_j+2)R}\right)^{\frac{q((pq)^j-1)}{pq-1}}
\end{array}
\]
because of  
\[
\dy\frac{t-\rho-R}{2}-\dy\frac{2a_j+2}{4a_j+5}(t-\rho-R)
=\dy\frac{1}{2(4a_j+5)}(t-\rho-R)\ge\dy\frac{t-\rho-R}{12a_j}. 
\]
Hence we obtain
\[
\begin{array}{ll}
U''(t)
&\ge\dy\frac{CC_j'^p}{(nq+n+1)^p(12a_j)^{npq+(n+1)p}}\\
&\quad\times
\dy\int_0^{t-(4a_j+6)R}\frac{\rho^{n-1-(n-1)p/2}(t-\rho-R)^{npq+(n+1)p}}
{(t-\rho +R)^{(n-1)p(pq+2q+1)/2}}d\rho\\
&\hspace{100pt}\times\left(\log\dy\frac{t-\rho+(8a_j^2+22a_j+14)R}
{2(4a_j+5)(a_j+2)R}\right)^{\frac{pq((pq)^j-1)}{pq-1}}\\
&\ge\dy\frac{CC_j'^p\{t-(4a_j+6)R\}^{-(n-1)p/2}}{(nq+n+q)^p(12a_j)^{npq+(n+1)p}}\\
&\quad\times
\dy\int_0^{t-(4a_j+6)R}\frac{\rho^{n-1}\{t-\rho-(4a_j+5)R\}^{npq+(n+1)p}}
{(t-\rho +R)^{(n-1)p(pq+2q+1)/2}}d\rho\\
&\hspace{100pt}\times
\left(\log\dy\frac{t-\rho+(8a_j^2+22a_j+14)R}
{2(4a_j+5)(a_j+2)R}\right)^{\frac{pq((pq)^j-1)}{pq-1}}.
\end{array}
\]
Since 
\[
t-\rho+R\le (4a_j+7)\{t-\rho-(4a_j+5)R\}\le 6a_j\{t-\rho-(4a_j+5)R\}
\]
is valid for $t-\rho\ge (4a_j+6)R$ and 
\[
\frac{n-1}{2}p(pq+2q+1)-npq-(n+1)p=1
\]
holds for $F(p,q,n)=0$, the $\rho$-integral is dominated from below by
\[
\begin{array}{l}
\dy\frac{\{t-(4a_j+6)R\}^{n-1}}{2^{n-1}(6a_j)^{(n-1)p(pq+2q+1)/2}}\\
\quad\times \dy\int_{(t-(4a_j+6)R)/2}^{t-(4a_j+6)R}
\dy\frac{d\rho}{t-\rho-(4a_j+5)R}
\left(\log\dy\frac{t-\rho+(8a_j^2+22a_j+14)R}
{2(4a_j+5)(a_j+2)R}\right)^{\frac{pq((pq)^j-1)}{pq-1}}\\
\ge\dy\frac{\{t-(4a_j+6)R\}^{n-1}}{2^{n-1}(6a_j)^{(n-1)p(pq+2q+1)/2}}\\
\quad\times \dy\int_{(t-(4a_j+6)R)/2}^{t-(4a_j+6)R}
\dy\frac{d\rho}{t-\rho+(8a_j^2+22a_j+14)R}
\left(\log\dy\frac{t-\rho+(8a_j^2+22a_j+14)R}
{2(4a_j+5)(a_j+2)R}\right)^{\frac{pq((pq)^j-1)}{pq-1}}\\
\ge\dy\frac{(pq-1)\{t-(4a_j+6)R\}^{n-1}}{2^{n-1}(6a_j)^{(n-1)p(pq+2q+1)/2}(pq)^{j+1}}
\left(\log\dy\frac{t+(16a_j^2+48a_j+34)R}
{4(4a_j+5)(a_j+2)R}\right)^{\frac{(pq)^{j+1}-1}{pq-1}}.
\end{array}
\]
Setting $a_{j+1}=4a_j+6$, we get the desired inequality for $j+1$;
\[
U''(t)\ge C_{j+1}(t-a_{j+1}R)^{n-1-(n-1)p/2}
\left(\log\dy\frac{t+(a_{j+1}^2-2)R}
{(a_{j+1}-1)(a_{j+1}+2)R}\right)^{\frac{(pq)^{j+1}-1}{pq-1}},
\]
where we set 
\[
C_{j+1}
=\frac{(pq-1)CC_j'^p}{(nq+n+1)^p2^{npq+(n+1)p+n-1}
(6a_j)^{(n-1)p^2q/2+(2n-1)pq+(3n+1)p/2}(pq)^{j+1}}.
\]
\par
To end the proof, we shall fix all the coefficients, $C_j$. 
It follows from $a_j=6\cdot 4^{j-1}-2\le 2\cdot 4^j$
and the definition of $C_j'$, (\ref{Cj'}), that 
\[
C_{j+1}= \dy\frac{MC_j^{pq}}{N^j},
\]
where 
\[
\begin{array}{l}
M=\dy\frac{(pq-1)C^{p+1}}{pqn^p(nq+n+1)^p2^{(n-1)p^2q+6npq+(5n+2)p+n-1}
3^{(n-1)p^2q/2+(3n-1)pq+(3n+1)p/2}},\\
N=4^{(n-1)p^2q/2+2npq+(3n+1)p/2}pq. 
\end{array}
\]
This equality is rewritten as 
\[
\log C_{j+1}=pq\log C_j+\log M-j\log N. 
\]
Then, one can easily get 
\[
\begin{array}{ll}
\log C_{j+1}
&=(pq)^j\log C_1+\dy\sum_{k=1}^j(pq)^{j-k}\log M-\dy\sum_{k=1}^jk(pq)^{j-k}\log N\\
&=(pq)^j(\log C_1+S_{j+1}), 
\end{array}
\]
where we set 
\begin{equation}
\label{S}
S_j=\dy\sum_{k=1}^{j-1}\frac{\log M-k\log N}{(pq)^k}. 
\end{equation}
Note that $S_j$ converges as $j\to\infty$. 
Therefore this completes the proof of Proposition \ref{prop:iteration1}. 
\par 
We omit to show the proof of Proposition \ref{prop:iteration2}
because it is almost the same as the single case in Takamura and Wakasa \cite{TW}.
The difference from the proof Proposition \ref{prop:iteration1} appears
in handling of logarithmic terms. 
In order to prove Proposition \ref{prop:iteration2}, 
we should integrate the logarithmic term at every steps in the iteration 
while such an integration is required only to get the estimate for $U''(t)$
in the proof of Proposition \ref{prop:iteration1}. 
%%%%%%%%%%%%%%%%%%%%%%%%%%%%%%%%%%%%%%%%%%%%%%%%

\par
\begin{center}
{\bf 
Acknowledgment}
\end{center}
\par
This manuscript was partially prepared during the first author's stay
at University of Pisa in Italy from 4/4/2011 to 28/9/2011
as an overseas research fellow sponsored
by Institute of National College of Technology, Japan. 
She is deeply grateful to Professor Vladimir Georgiev for his hearty hospitality,
a lot of help on her stay and many fruitful discussions. 
She also thanks to all the members of Department of Mathematics,
University of Pisa for preparations of the necessaries for her activities. 
%%%%%%%%%%%%%%%%%%%%%%%%%%%%%%%%%%%%%%%%%%%%%%%

\end{document}